\documentclass[12pt,preprint]{elsarticle}

\usepackage{amssymb}

\usepackage[utf8]{inputenc}
\usepackage{graphicx}
\usepackage{dcolumn}
\usepackage{dsfont}
\usepackage{amsthm}
\usepackage{epstopdf}
\usepackage{mathtools}
\usepackage{textcomp}
\usepackage{colortbl}
\usepackage{float}
\usepackage[]{algorithm2e}
\usepackage{bm}
\usepackage[format=plain,font={small},justification=centerlast]{caption}
\usepackage{subcaption}
\usepackage{color}
\usepackage{tabularx,ragged2e,booktabs}
\usepackage{MnSymbol}
\usepackage{wasysym}
\usepackage{anysize}
\usepackage{epstopdf}
\usepackage{verbatim}
\usepackage{tikz}
\usepackage{pgfplots}
\usepackage{mathtools}
\usepackage{natbib}
\usepackage{xcolor}
\DeclareGraphicsExtensions{.pdf,.png,.jpg,.pdf}



\biboptions{comma,square}

\newcommand{\bc}{\begin{center}}
	\newcommand{\ec}{\end{center}}
\newcommand{\be}{\begin{equation}}
\newcommand{\ee}{\end{equation}}
\newcommand{\bea}{\begin{eqnarray}}
\newcommand{\eea}{\end{eqnarray}}
\newcommand{\ba}{\begin{array}}
	\newcommand{\ea}{\end{array}}

\newcommand{\x}{\mathbf{x}}

\newcommand{\m}{\mathbf{m}}
\newcommand{\dd}{\mathbf{d}}
\newcommand{\dbf}{\mathbf{d}}

\pgfplotsset{compat=1.14}
\journal{Elsevier}

\begin{document}
	
	\begin{frontmatter}
		
		
		\title{Application of an RBF-FD solver for the Helmholtz equation to full-waveform inversion}


		
		\author[1,2]{Mauricio A. Londo\~{n}o}
		\ead{alejandro.londono@udea.edu.co, malondono@unal.edu.co}
		
		\author[2]{Francisco J. Rodríguez-Cortés}
		\ead{frrodriguezc@unal.edu.co}
		
		\address[1]{Instituto de Matemáticas \\
			Universidad de Antioquia}
			
		\address[2]{Escuela de Estadística \\
			Universidad Nacional de Colombia, sede Medellín}
		
		\begin{abstract}
	Full waveform inversion (FWI) is one of a family of methods that allows the reconstruction
	of earth subsurface parameters from measurements of waves at or near the surface. This
	is a numerical optimization problem that uses the whole waveform information of all
	arrivals to update the subsurface parameters that govern seismic wave propagation. We
	apply FWI in the multi-scale approach on two well-known benchmarks: Marmousi and
	2004 BP velocity model. For the forward modeling, we use an RBF-FD solver on hexagonal
	grids and quasi-optimal shape parameters, developed in \cite{londono2019}.
		\end{abstract}
		
		
	\end{frontmatter}
	
	
	\section{Introduction}
    
    One of the most important challenges in seismic exploration is to get an accurate model of the wave propagation velocity of subsurface earth. Currently, the Full Waveform Inversion method (FWI) is a tool with great acceptance and constant evolution. Details can be consulted in \cite{introducttion_to_FWI_vireux_2017}.


In a 2D acoustic seismic acquisition are placed, near to ground surface, a set of sources $\{\x_{s_i}\}$ (\emph{shots}) and a set of receivers $\{\x_{r_j}\}$ (\emph{geophones or hydrophones}). Data is recorded through the receivers (\emph{observed data}\footnote{Actually, data are recorded in time domain then it could be transformed from time to frequency domain by means of Fourier transform.}) forming  arrays  $\mathbf{d_{\omega}}=\{d^{s_{i},r_{j}}_{\omega}\}$, for certain frequencies $\omega=\omega_0,\ \omega_1,\ \ldots,\ \omega_l$; where $\omega_0<\omega_1<\cdots<\omega_l$. It is supposed that data $d^{s_{i},r_{j}}_{\omega}=u^{s_i}_{\omega}(\x_{r_j})$, where  $u_{\omega}^{s_i}$ (\emph{the pressure wave-field}) satisfy the Helmholtz equation
\begin{equation*}
-\Delta u^{s_i}_{\omega}(\x)-\omega^2c(\x)^{-2}u^{s_i}_{\omega}(\x)=\delta(\x-\x_{s_i}).
\end{equation*}
 Data  are the response produced by subsurface geologic structures when the acoustic medium, with unknown velocity $c(\x)$ (\emph{true model}), is perturbed with a certain source energy (a shot). In what follows we take $\m=c(\x)^{-2}$ by notational convenience suggested for the literature about FWI.  A typical FWI scheme consists in finding a good approximation of $\m$ from known data $\dd_{\omega}$, by  minimizing a certain misfit functional, for example, we use the classic least squares
\begin{equation*}
\mathcal{J}_{\omega}(\m)=\frac{1}{2}\|\dbf_{\omega}-\mathcal{F}_{\omega}(\m)\|^2_2,
\end{equation*}
where $\mathcal{F}_{\omega}$, called the \emph{forward modeling map}, is defined by the array $\mathcal{F}_{\omega}(\m)=\{u_{\omega}^{s_{i}}(\x_r)\}$ and 
\begin{equation*}
-\Delta u^{s_i}_{\omega}(\x)-\omega^2\m\, u^{s_i}_{\omega}(\x)=\delta(\x-\x_{s_i}).
\end{equation*}

 Note that to compute $\mathcal{F}_{\omega}(\m )$ is necessary to solve a Helmholtz problem with many single sources.
The frequency domain acoustic FWI method can be viewed as the optimization problem:
\begin{equation}\label{eq:optimization_problem}
\mbox{Given}\ \ \dd_{\omega},\  \mbox{ to find } \arg\min_{\m} \mathcal{J}_{\omega}(\m).
\end{equation}
Typical solutions for this problem are given by iterative methods such as
the gradient descent method
\begin{equation}\label{eq:gradiend_descent}
\m^{(k+1)}=\m^{(k)}-\alpha_k\nabla\mathcal{J}_{\omega}(\m^{(k)})
\end{equation}
starting with a "good" initial model $\m^{(0)}$ and where $\alpha_k$ is the step length, which can be computed efficiently by using either Barzilai-Borwein method (BB) \cite{BB_Fletcher2005}, BFGS  or L-BFGS methods \cite{optimization_nocedal_2006}. To compute $\nabla\mathcal{J}_{\omega}(\m^{(k)})$ we have used the adjoint operator  of the forward modeling operator $\mathcal{F}$. In geophysics literature this is known as the adjoint state method, details about it can be consulted in \cite{Plessix2006_A_review_of_the_adjoint_state_method}.

Since the the objective function is not convex , in general, we are faced with multimodal extremes. Hence in solving the FWI problem  we can  be trapped in local minima points which could take us to a wrong solution. Success of FWI depends strongly on the initial model,  which must offer a response,  in low frequency, close to the real model.
 
The availability of low-frequency data  is an important factor to the success of FWI. The low frequencies help to determine the kinematically relevant low-wavenumber components of the velocity model, which are in turn needed to avoid convergence of FWI to spurious local minima. However, acquiring data less than 2Hz or 3Hz from the acquisition of real data  is a  challenging and expensive task. Some works have explored the possibility of synthesizing the low frequencies computationally from high-frequency data and use the resulting prediction of the missing data to seed the frequency sweep of FWI for the multi-scale approach \cite{li-demanet}, overcoming the cycle skipping phenomena \cite{2018_Full-waveform_inversion_using_a_nonlinearly_smoothed_wavefield}, \cite{2017_low_frequency_envelope_FWI}.
 
  Normally, due to the lack of low-frequency components, algorithms for FWI are tested with initial models from either smoothed or filtered versions of the real model. However, if low frequency data are available it can be used simple initial such as a linear one, reducing the risk of getting trapped at a local minimum. 
  

We assume prior information about true model $\m$, such as its minimum and maximum values, and shallow information  which in practice can be known from prior geological studies. We start with a depth linear initial model $\m_0$ created from near values to minimum and maximum values of true model $\m$. The initial frequency $\omega_0$ for multi-scale method, which in our tests have been successful, is approximated by $$\omega_0\approx\frac{2\pi}{z_d\sqrt{\min(\m_0)}},$$
 where $z_d$ is the depth of the target model  and $[\m_0]$ is the mean value of $\m_0$. 
We use a multi-scale approach with single frequencies, i.e., for the frequency $\omega_0$ and the initial model $\m_0$, we obtain a model $\m_{\omega_0}$ by solving \eqref{eq:optimization_problem} using the descent method $\eqref{eq:gradiend_descent}$ with the observed data $\dd_{\omega_0}$. The found model, $\m_{\omega_0}$, is used as initial model to invert data $\dd_{\omega_1}$. In genereal, if $\m_{\omega_j}$ is the accepted final model using the data $\dd_{\omega_j}$, then $\m_{\omega_j}$ is used as initial model to invert data $\dd_{\omega_{j+1}}$. This process is summarized in the algorithm \ref{alg:FWI_multi-scale}.

\RestyleAlgo{boxruled}
\begin{center}	
\begin{algorithm}[h]
	\small
	\KwData{An initial model $\m_0$, and for some selected frequencies, observed data: $\dd_{\omega_0},\dd_{\omega_1},\ldots,\dd_{\omega_l}$}
	\KwResult{An opproximation of the true model: $\m\approx\m_{\omega_{l+1}}$ }
	initialize $\m_{\omega_0}=\m_0$\; 
	\For{p=0,\ldots,l}{
		define: $tol_g>0$, $tol_J>0$, $maxiter>0$\;
		initialize:	$k=0$\;
		$\m^{(k)}=\m_{\omega_p}$\;
		$\omega=\omega_p$\;	
	\While{$k<$  maxiter \emph{ and} $\|\nabla\mathcal{J}_{\omega}(\m^{(k)})\|>tol_g$ \emph{and} $\mathcal{J}_{\omega}(\m^{(k)})>tol_J$ }{
		 compute $\alpha_k$ with some optimizer (e.g. BB, BFGS, L-BFGS, etc)\;
		$\m^{(k+1)}=\m^{(k)}-\alpha_k\nabla\mathcal{J}_{\omega}(\m^{(k)})$\;
		$k=k+1$\;		
}
$\m_{\omega_{p+1}}=\m_{k+1}$\;
}\normalsize
\caption{Multi-scale FWI by single frequency.\vspace{1cm}}
\label{alg:FWI_multi-scale}
\end{algorithm}
\end{center}
\normalsize

We use RBF-FD on uniform hexagonal grid with mesh size adapted to the minimum wavelength  to save computational cost at low frequencies \cite{londono2019}. In this way it is feasible to consider a large wavelength to obtain adequate initial models. 
\subsection{Barzilai-Borwein optimization method}

To compute the step length $\alpha_k$ in \eqref{eq:gradiend_descent} we have chosen the Barzilai-Borwein (BB) optimization method \cite{BB_Fletcher2005}. BB is a steepest descent method based on non-monotone line search technique. BB is suitable for large scale problems. One important feature about BB is that it does not require the Hessian matrix to get $\alpha_k$ and this make it fast and relatively cheap for each iteration. However, the sequence of error values $\mathcal{J}(\m^{(k)})$ is not monotone. \cite{FWI_BB} shown evidence that BB is effective in order to save memory for Frequency domain FWI, with good trade-off between cost and benefit.

\subsection{Limited-Memory BFGS}
Within quasi-Newton algorithms one of the most popular is the BFGS method and its variant Limited-Memory BFGS (L-BFGS), which for FWI is suitable because the steep length is calculated by approximations of the Hessian through evaluations of the gradient. We have used the matlab library by \cite{page:lbfgs}. Details can be seen in \cite{optimization_nocedal_2006}.

\subsection{Numerical Tests } 
In this part we describe initial settings to perform FWI on --offshore models-- Marmousi and 2004 BP. For the forward modeling we use GRBF-FD7p for PML \cite{londono2019}. In all tests we produce synthetic observed data with the  GRBF-FD7p solver.

\subsubsection{Offshore models}

\begin{figure}[H]
	\begin{center}
		\begin{tabular}{l}		
			\includegraphics[width=8.1cm]{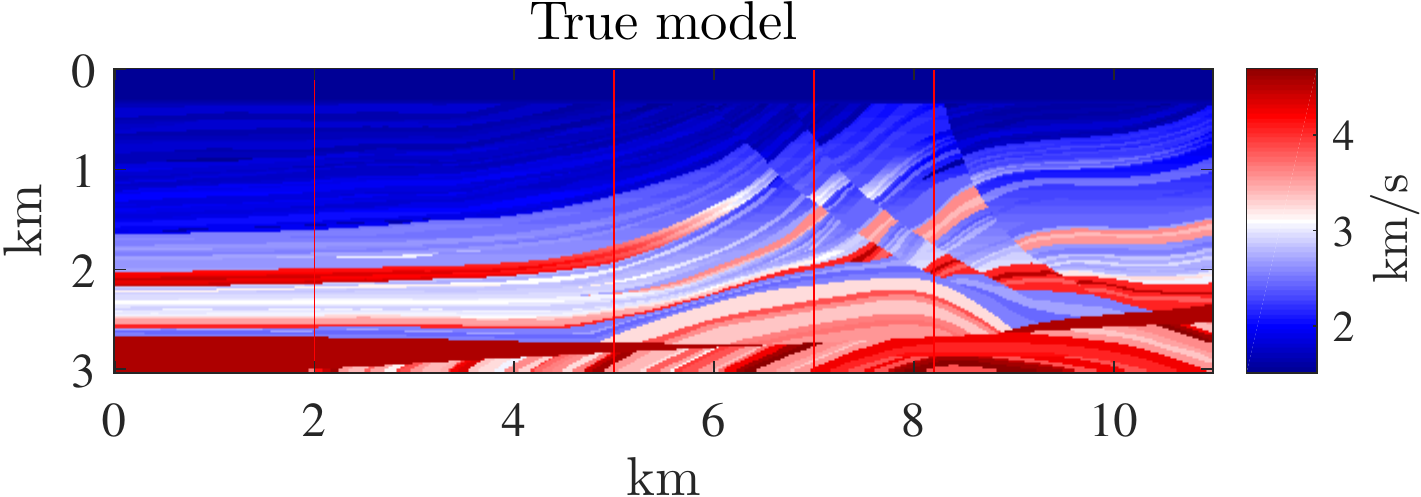}\vspace{-5mm}\\ \includegraphics[width={\textwidth}]{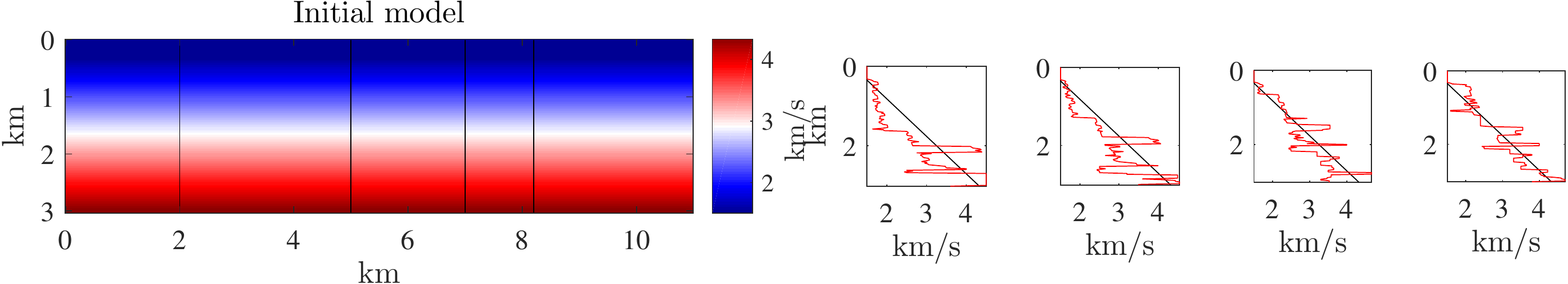}\vspace{-5mm} \\		
			\includegraphics[width={\textwidth}]{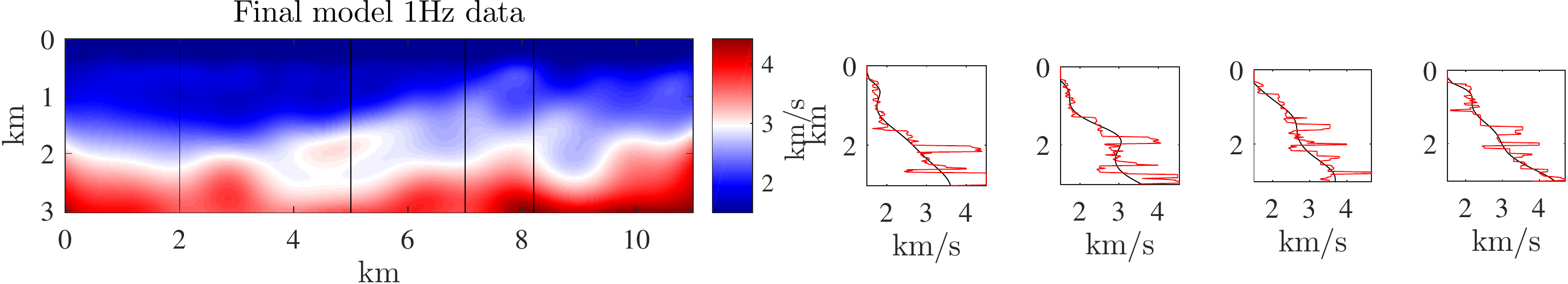}\vspace{-5mm} \\ \includegraphics[width={\textwidth}]{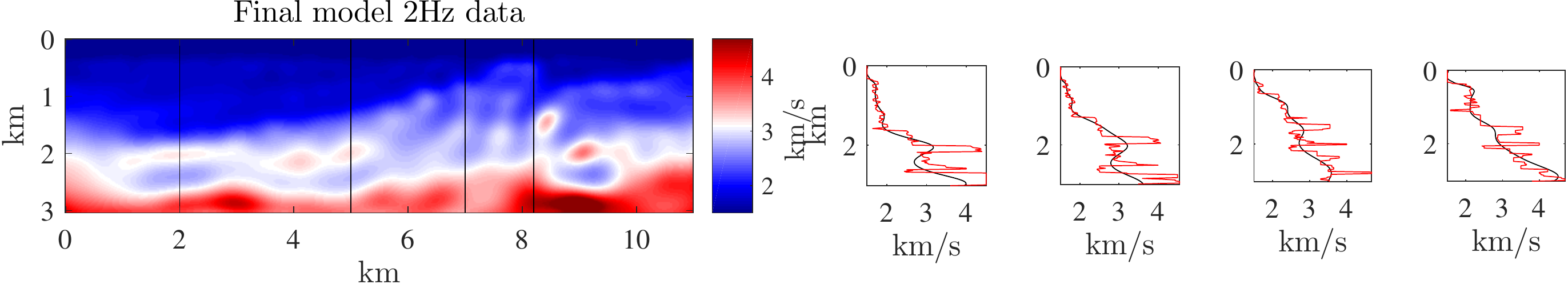}\vspace{-5mm} \\
			\includegraphics[width={\textwidth}]{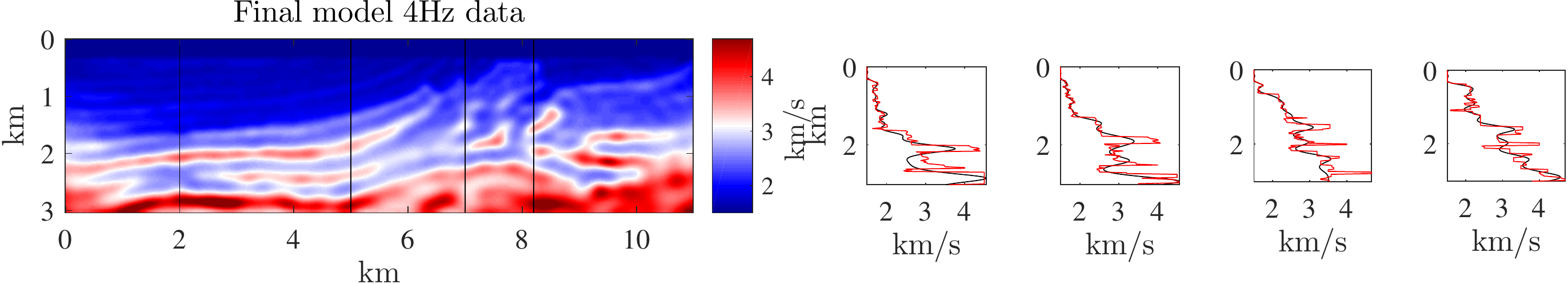}\vspace{-5mm} \\
			\includegraphics[width={\textwidth}]{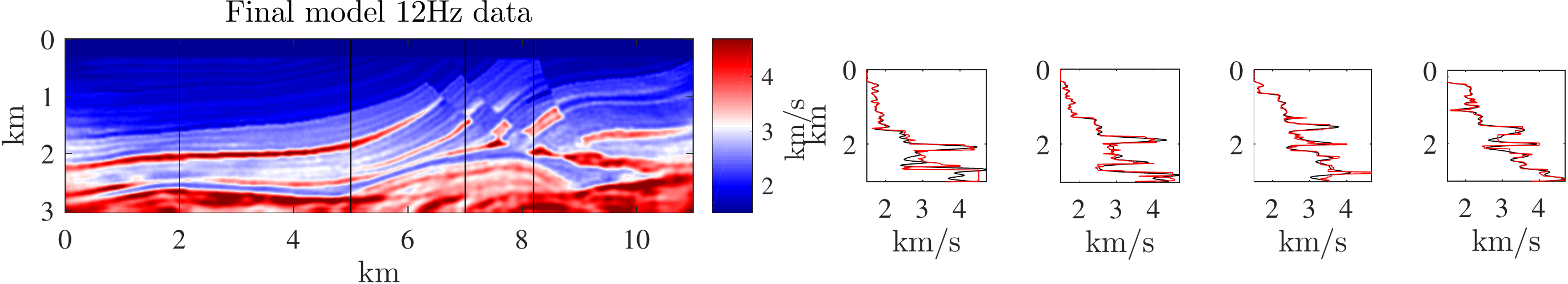}\vspace{-5mm} \\
			\includegraphics[width={\textwidth}]{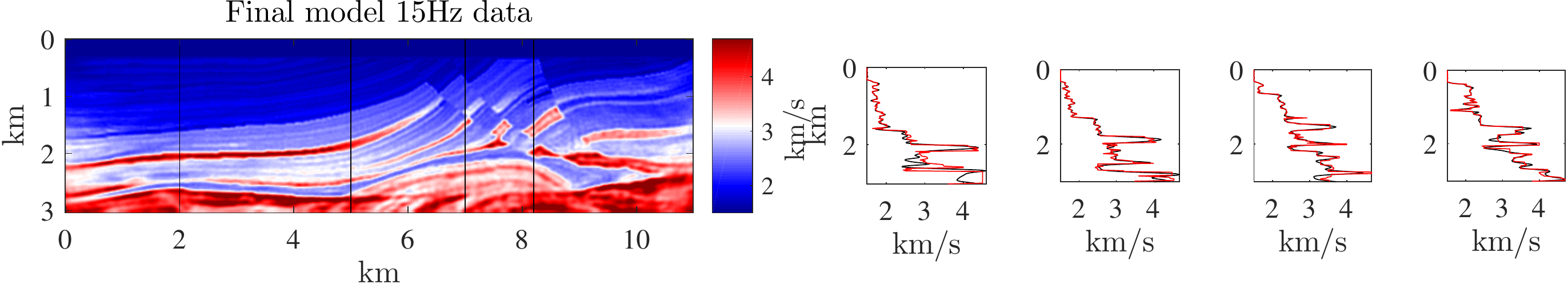} 
		\end{tabular} \caption{Marmousi model: Some resulting models by using the multi-scale approach with 15 single frequencies from 1Hz to 15Hz, with step frequency of 1Hz. The group of plots on the right corresponds to a comparison between vertical profiles of true model (red lines) and resulting  models (black lines), taken at $x=2 \mbox{km},\, 5 \mbox{km},\, 7\mbox{km},\,\mbox{and } 8.2\mbox{km}$.}
		\label{fig:marm_FWI}
	\end{center}
\end{figure}

\textbf{Marmousi:}
We use a version of the Marmousi model \cite{marmousi_paper} given in a uniform square grid, in a matrix of size $n_z\times n_x=150\times550$ with horizontal and vertical separation of 20m. Thus, this model represents 3km in depth and 11km in surface. For this test we use 55 sources $\x_{s_i}$ uniformly distributed with a separation of 200m, where the first one is located to 20m from the edge. We used 276 receivers $\x_{r_j}$ uniformly distributed with a separation of 40m, where the first one is located on the edge, i.e., at $n_x=0$. We select frequencies $f=\frac{\omega}{2\pi}=1\mbox{Hz},2\mbox{Hz},\ldots,15\mbox{Hz}$. At each frequency the uniform hexagonal grid is built with $Ng=8.5\,$NPW for the minimum wavelength of the model, that is $\displaystyle\lambda_{min}=\min(c)/f$.  For the initial model we assume that there is known information of the first layer, which would correspond to seawater. We have used the  L-BFGS optimizer. Results of this test can be seen in Fig. \ref{fig:marm_FWI}.
	
\begin{figure}[H]
	\begin{center}
		\begin{tabular}{l}		
			\includegraphics[width=10.3cm]{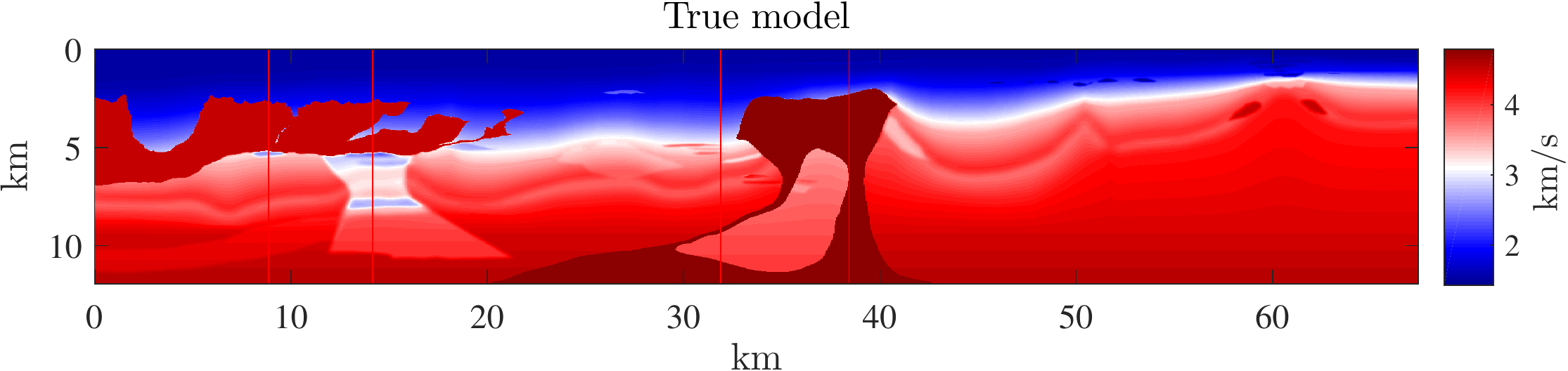}\vspace{-4mm}\\ \includegraphics[width=10.3cm]{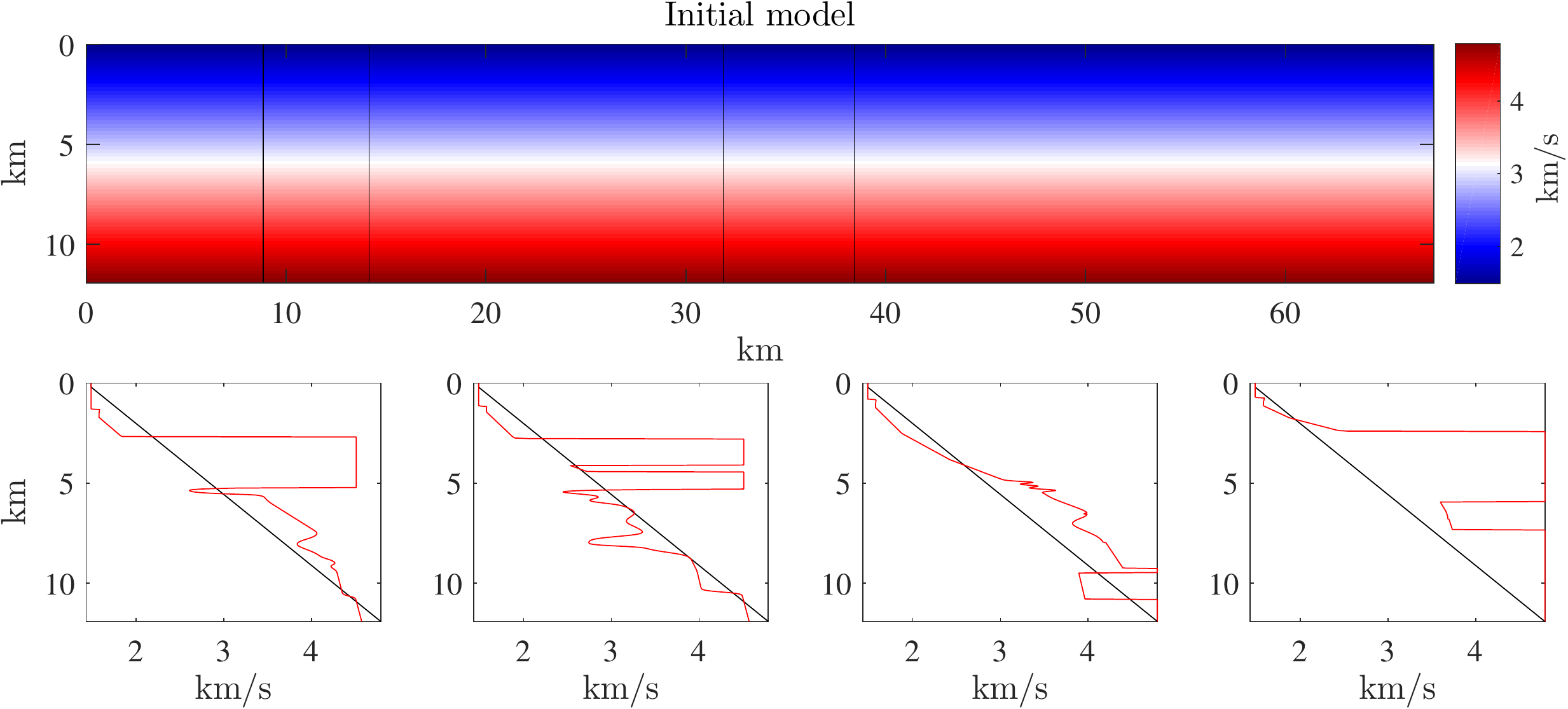}\vspace{-5mm}\\
			\includegraphics[width=10.3cm]{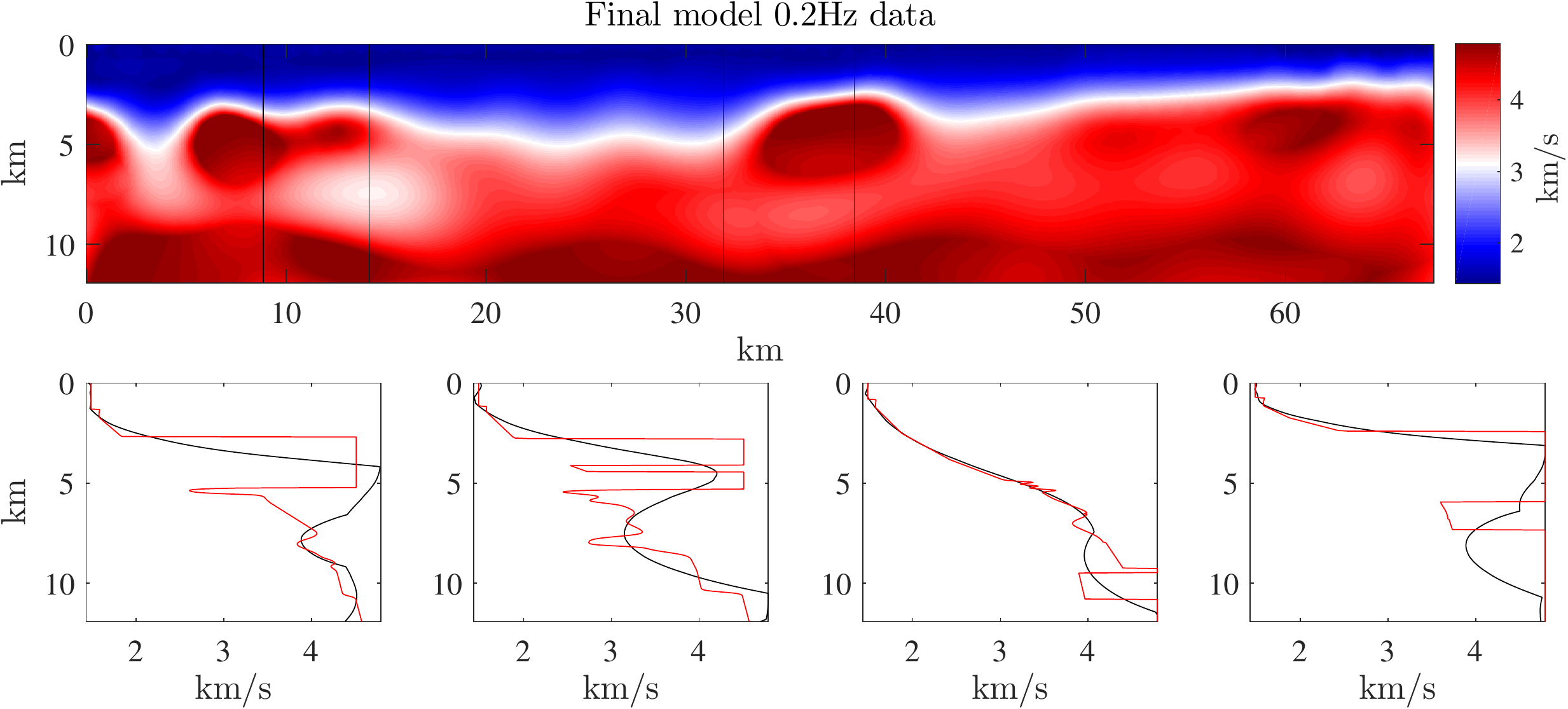} \vspace{-5mm}\\ \includegraphics[width=10.3cm]{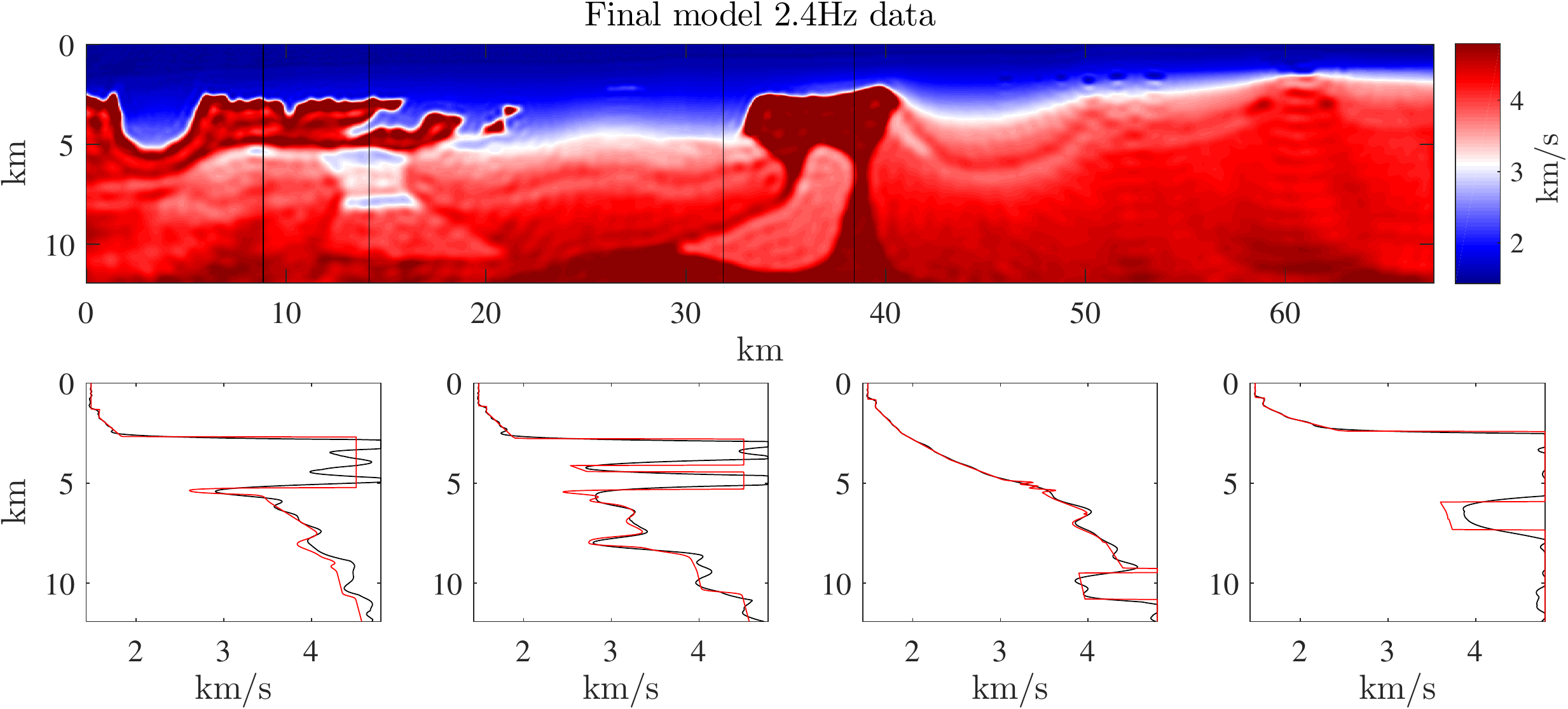}\vspace{-5mm} \\
			\includegraphics[width=10.3cm]{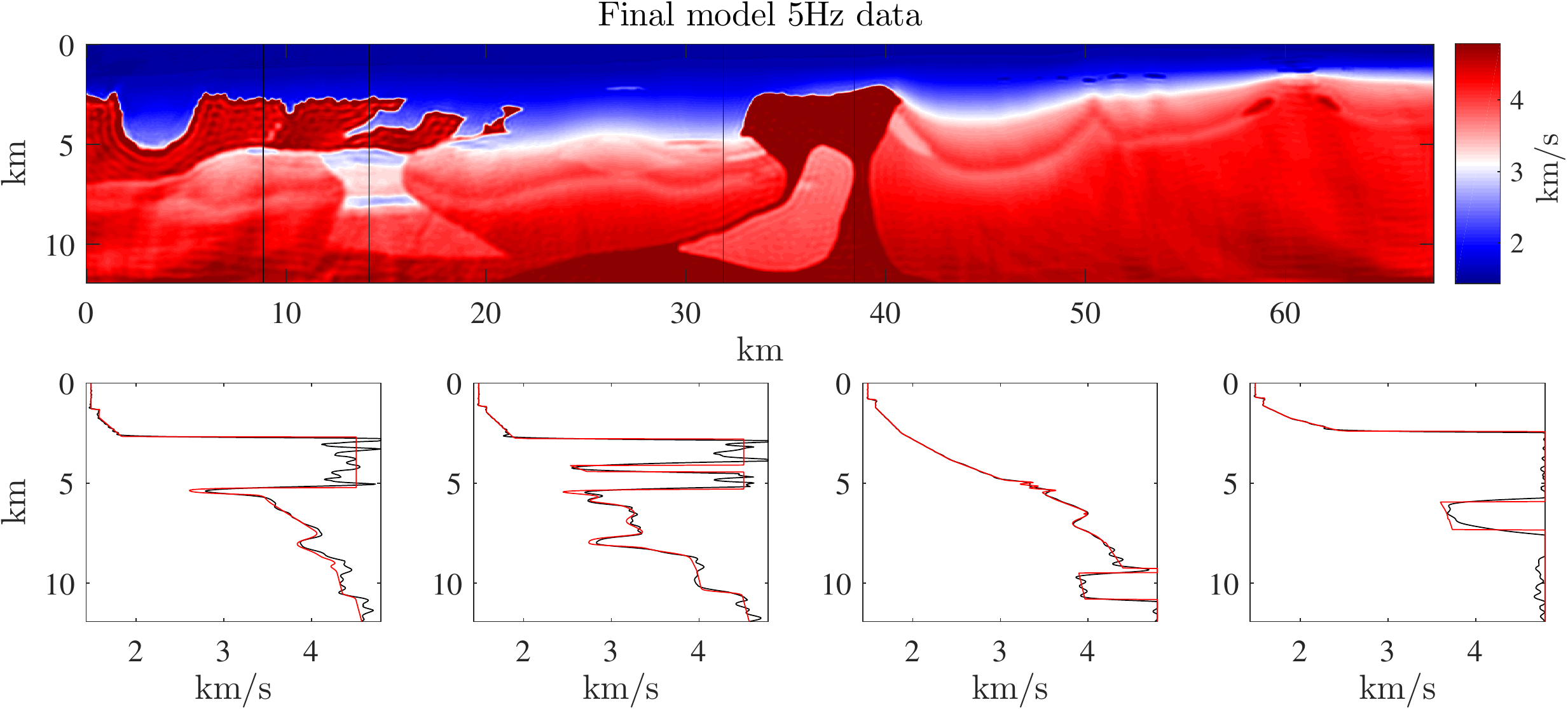}
		\end{tabular} \caption{2005 BP velocity model: Some resulting models using single frequency multi-scale approach with 25 frequencies from 0.2Hz to 5Hz. The group of plots below each model corresponds to  comparisons between vertical profiles of the true model (red lines) and resulting  models (black lines), taken at $x=8.85 \mbox{km},\, 14.16 \mbox{km},\, 31.88\mbox{km},\,\mbox{and } 38.42 \mbox{km}$. To minimize the misfit function we used the BB optimizer. Results are shown in Table \ref{table:2004bp_inversion}.}
		\label{fig:bp2004_inversion}
	\end{center}
\end{figure}

\noindent\textbf{2004 BP:} For this test we used a reduced version of 2004 BP given in a uniform square grid of $n_z\times n_x=478\times2698$, with separation of 25m which represents 11.925km in depth and 67.425km in surface. We used 270 sources $\x_{s_i}$ uniformly distributed with separation of 250m, where the first one is located to 25m from the edge. We used 1350 receivers $\x_{r_j}$ with separation of 50m, where the first one is located on the edge, i.e., at $n_x=0$. We have selected frequencies $f=\frac{\omega}{2\pi}=0.2\mbox{Hz},\ 0.4\mbox{Hz},\ 0.6\mbox{Hz},\ \ldots,\ 5.0\mbox{Hz}$; with frequency step of $0.2\mbox{Hz}$. At each frequency the uniform hexagonal grid is built with $N_g=6.5\, $NPW, in according to the minimum wavelength $\displaystyle\lambda_{min}=\min(c)/f$. For the PML thickness $\delta$,  we take $\delta=\lambda=[c]/f$, where  $[c]$ is the mean value of $c$. For the initial model we assume that there is known information of a shallow layer, which would correspond to sea water. Results for this inversion test can be seen in Fig. \ref{fig:bp2004_inversion} and Table \ref{table:2004bp_inversion}.

\begin{table}[H]
	\centering
	\begin{tabular}{ccccc}
		\hline\hline
		$\omega/2\pi$ (Hz) & Time per iteration (s) &\# iterations & Size of $\mathbf{H}$ (inner nodes) & $\|\nabla \mathcal{J}_{\omega}(\m)\|$ \\
		\hline 
		\hline
		$0.2$ & 3.40 & 748 & 1985 & 1.99e-8 \\ 
		$0.4$ & 6.42 & 116 & 5285 & 1.29e-7 \\ 
		$0.6$ & 9.42 & 279 &10022 & 1.39e-7 \\ 
		$0.8$ & 12.8 & 237 &16738 & 1.31e-7 \\
		$1.0$ & 17.3  & 146  &24602& 1.38e-7 \\
		$1.2$ & 25.2 & 64  & 33865&  3.58e-7  \\
		$1.4$ & 33.4  & 85 & 44592 & 3.93e-7 \\
		$1.6$ & 43.8  & 40 & 57903 & 3.69e-7  \\
		$1.8$ & 57.4  & 46 & 71692 & 3.54e-7 \\
		$2.0$ & 64.6 & 44 & 87084 & 3.62e-7 \\ 
		$2.2$ & 80.3  & 50  & 103813 &  3.87e-7  \\ 
		$2.4$ & 98.8  & 40 & 123503 & 3.98e-7 \\ 
		$2.6$ & 120  & 39  & 143282 & 4.00e-7  \\ 
		$2.8$ & 132  & 56 & 164714 & 3.81e-7  \\ 
		$3.0$ & 158 & 47 & 189501 & 3.93e-7 \\ 
		$3.2$ & 182 & 39 & 213818& 3.31e-7  \\ 
		$3.4$ & 214  & 29  & 239599 & 3.91e-7  \\ 
		$3.6$ & 242  & 42  & 266844 & 3.93e-7  \\ 
		$3.8$ & 266  & 51  & 298161  & 3.66e-7  \\ 
		$4.0$ & 291  & 54  & 328468  & 3.35e-7   \\ 
		$4.2$ & 319 & 30  & 360239  & 3.72e-7  \\ 
		$4.4$ & 355 & 40  & 393474  & 3.05e-7 \\ 
		$4.6$ & 397  & 53  & 431321 & 3.76e-7  \\ 
		$4.8$ & 438  & 32  & 467618 & 3.57e-7  \\ 
		$5.0$ & 490  & 35  & 505706  & 3.80e-7   \\ 
		\hline\hline
	\end{tabular}%
	\caption{Results for FWI applied to 2004 BP velocity model.}	
	\label{table:2004bp_inversion}
\end{table}

\section{Conclusions}

With Gaussian RBF we have exploited the prior known of the local wavenumber $k=\omega c^{-1}$ to obtain near-optimal shape parameters to approximate the Laplace operator in accord to wavelength. Specifically, for nodes distributed in a hexagonal grid, we have developed an scheme with closed formulas of near-optimal weights. In this approach the, GRBF-FD7p is a efficient method which shows a good trade-off between cost and benefit, which allowed to apply GRBF-FD to the frequency domain FWI method obtaining satisfactory results.

	\section*{Acknowledgments}
	
This work is supported by MINCIENCIAS-Colombia as a part of the research project grant No. 80740-735-2020.

	\bibliographystyle{apalike}
	\bibliography{references.bib}

	
	
	
	
	

\end{document}